\documentclass[12pt]{article}
\usepackage{custarticle}

\title{New examples of $\Z_2$-harmonic 1-forms and their deformations}
\author{Andriy Haydys \\ Université Libre de Bruxelles
	\and Rafe Mazzeo \\ Stanford University
	\and Ryosuke Takahashi \\ National Cheng Kung University}
\date{10 FEB 2025}

\begin{document}
	\maketitle

\begin{abstract}
We collect a number of elementary constructions of $\Z_2$ harmonic $1$-forms, and of families of these objects. 
These examples show that the branching set $\Sigma$ of a $\Z_2$ harmonic 1-form may exhibit the following features: 
i) $\Sigma$ may be a non-trivial link; ii) $\Sigma$ may be a multiple cover; iii) $\Sigma$ may be immersed, and appear as a 
limit of smoothly embedded branching loci; iv) there are families of $\Z_2$ harmonic $1$-forms whose branching sets $\Sigma$
have tangent cones filling out a positive dimensional space, even modulo isometries. 
We show that Features i) and ii) occur already in dimension three, while the remaining ones appear at least in dimension four and higher. 
\end{abstract}
	
\section{Introduction}
Let $M$ be an oriented Riemannian $n$-manifold. Any real Euclidean line bundle $\cI\to M$ is equipped with a flat metric 
connection $\nabla = \nabla^{\cI}$ uniquely specified by the requirement $\nabla s =0$ for each smooth local section of 
constant pointwise norm.  We consider here sections of this bundle, or $1$-forms with values in $\cI$, i.e., sections
of $T^*M \otimes \cI$. If this flat connection has nontrivial monodromy, then these sections are simply ordinary
functions or $1$-forms in a local trivialization, but are well-defined only up to a factor of $\pm 1$. We shall be interested
in the the situation where $\cI$ is defined only over an open set $M \setminus \Sigma$, where $\Sigma$ is a closed subset.
Of particular interest is the case where $\Sigma$ is of codimension $2$, which makes it possible for $\cI$ to have nontrivial monodromy.

In terms of this connection and the metric $g$ on $M$ we may define harmonic sections of $\cI$ and harmonic $1$-forms 
with values in $\cI$. These may be described, alternately, as harmonic functions or $1$-forms defined on $M \setminus \Sigma$
which are defined up to a sign $\pm 1$. Somewhat more generally, we may consider harmonic sections of other bundles over $M$
coupled to $\cI$; of particular interest are harmonic sections of $\slS \otimes \cI$, where $\slS$ is the spin bundle. These are
called $\Z_2$ harmonic spinors.  

In this paper we focus mostly on $\Z_2$ harmonic functions and $1$-forms; we regard each such object as a triple $(u, \Sigma, \cI)$ 
or $(\om, \Sigma, \cI)$, where $\Sigma$ is a closed subset of $M$, $\cI$ is a real Euclidean line bundle over $M\setminus \Sigma$, 
and $u$ and $\om$ are sections of $\cI$ and $T^*M \otimes \cI$, respectively, such that 
\begin{equation}
\label{Eq_Z2Harm1Form}
d^* du = 0, \quad \mbox{and} \quad (d + d^* )\, \om= 0
\end{equation}
hold everywhere on $M\setminus \Sigma$. 
Here the operators $d$ and $d^*$ are the exterior derivative 
and its adjoint twisted by $\nabla^\cI$.
In the case of harmonic functions, we require also that both $|u|$ and $|du|$  extend to all of $M$ as H\"older-continuous
functions which vanish along $\Sigma$; in the case of harmonic 1-forms, only the H\"older-continuity of $|\om|$ makes sense.  
 In precisely the same way we consider $\Z_2$ holomorphic sections of $\cI \otimes \C$.
We call $\Sigma$ the branching set in all these cases. 

Although this definition makes sense in any dimension, we concentrate here on dimensions $n=3$ and $4$ due to the emerging
importance of $\Z_2$ harmonic functions and $1$-forms in these dimensions in various applications. 

The notion of $\Z_2$ harmonic functions, $1$-forms and spinors, was introduced and emphasized initially by Taubes~\cite{Taubes13_PSL2Ccmpt,Taubes15_CorrigendumToPSL2C,Taubes14_ZeroLoci_Arx}, although of course various
special cases are classical. Such objects are closely related to diverging sequences of flat stable $\PSL(2,\C)$-connections 
when $\dim M =3$. It was shown later that $\Z_2$ harmonic $1$-forms/spinors appear in a number of other gauge-theoretic problems, 
for example the Seiberg--Witten equations with multiple spinors both in dimensions $3$ and $4$~\cite{HaydysWalpuski15_CompThm_GAFA,Taubes16_SWDim4_Arx}, the complex anti-self-duality
equations~\cite{Taubes13_CxASD_Arx}, the Vafa--Witten equations~\cite{Taubes17_VafaWitten_Arx}, the 
Kapustin--Witten equations \cite{Taubes18_SequencesKapustWitten_Arx}, etc. 
Apart from gauge theory, $\Z_2$ harmonic forms/maps appear for example also in relation to branched coverings 
of special Lagrangian submanifolds~\cite{He23_BranchedDeform} and certain fibrations of $\rG_2$
manifolds~\cite{Donaldson17_AdiabLimitsG2}.  We mention also that multi-valued functions satisfying linear and
nonlinear elliptic equations have been studied for several decades by Almgren, and more recently De Lellis and others, for 
their role in the study of minimal submanifolds of codimension greater than $1$. 

At present, relatively little is known about constraints on the structure of branching sets.  There is a Fredholm 
deformation theory for $\Z_2$ harmonic spinors for which the branching set is a smooth curve 
\cite{Takahashi15_Z2HarmSpinors_Arx, Parker23_DeformZ2HarmSpinors_Arx}; the branching sets appearing
in these deformation families remain smooth. It was proved by Zhang~\cite{Zhang22_Rectifiability} that any branching
set $\Sigma$ is necessarily $(n-2)$-rectifiable.  However, the first examples where $\Sigma$ is not smooth (at least
in dimension $3$) appeared only quite recently~\cite{TaubesWu20_Examples}.  In another direction, it is shown 
in \cite{Haydys22_SWPSL2R} that if $\Sigma$ is a link and $M$ is a rational homology three-sphere, 
then $\Delta_\Sigma(-1)=0$, where $\Delta_\Sigma$ is the Alexander polynomial of $\Sigma$. However, until now no examples 
of $\Z_2$ harmonic forms where $\Sigma$ is a non-trivial link in a $3$-manifold have been constructed. This lack of examples 
is not only aesthetically dissatisfying, but in fact an obstruction to developing a more systematic study 
of $\Z_2$ harmonic functions/$1$-forms/spinors. 

The main and rather modest aim of this manuscript is to provide a collection of new explicit examples and methods of 
construction of $\Z_2$ harmonic $1$-forms, using essentially entirely elementary methods.  Although $\Z_2$ harmonic
$1$-forms are typically somewhat more interesting, any such $1$-form $\om$ can be written (at least locally) as
$du$ where $u$ is a $\Z_2$ harmonic function. (The usual multi-valuedness arising from periods of $\om$ around $\Sigma$
does not arise because we insist that $|\om|$ vanish on $\Sigma$, though of course $u = \int \om$ exhibits the 
same $\pm 1$ indeterminacy.) In any case, because of this we primarily discuss $\Z_2$ harmonic functions. Note, however, 
that for $|\om|$ to be H\"older continuous across $\Sigma$, it is necessary that $|u|$ vanish to order greater than $1$
along this branching set, which is an added constraint.

In dimension $4$ we can take advantage of complex analytic techniques, and this leads easily to a wide variety of examples
which exhibit various types of interesting behavior. In particular, we show that the zero set of any holomorphic function
on a K\"ahler surface is the branching set of a $\Z_2$ harmonic $1$-form, and consequently, there exist one-parameter families 
$(\om_t,\Sigma_t, \cI_t)$ of $\Z_2$ harmonic $1$-forms such that $\Sigma_t$ is smoothly embedded when $t\neq 0$ but 
only immersed (or considerably more singular) when $t=0$. From this correspondence we can also see that the set of 
possible tangent cones at a singular point of $\Sigma$ has positive dimension, see Example~\ref{Ex_TangCones} for further details.

It is considerably harder to construct examples in dimension $3$, but using the theory of harmonic morphisms~\cite{BairdWood03_HarmMorphisms}, 
we construct (local) examples of $\Z_2$ harmonic $1$-forms for which the branching set is a Hopf link or any torus knot in $S^3$. We also 
find a family $(\om_a, \Sigma_a,\cI_a)$ such that $\Sigma_a$ is a knot, when $a \neq 0$, which is a non-trivial 
multiple cover of  the unknot $\Sigma_0$. 
This means that there exists a tubular neighborhood $U$ of $\Sigma_0$ such that the restriction of the projection $\pi\colon U\to \Sigma_0$ yields a covering $\Sigma_a\to \Sigma_0$ of finite degree $q$ for all  sufficiently small $a\neq 0$. In particular, $\Sigma_a$ converges to $q\Sigma_0$ as $a$ converges to $0$ (in the sense of currents, for example, after choosing suitable orientations). The reader may find more details in Example~\ref{Ex_TorusKnots}.  

\medskip

\textsc{Acknowledgments.} The first named author was partially supported by the ARC 
grant ``Transversality and reducible solutions in the Seiberg--Witten theory with multiple spinors'' of the ULB. The third named author was partially supported by Ministry of
Education, Taiwan, R.O.C. under Grant no. MOST 111-2636-M-006-023.
  
\section{Examples of $\Z_2$ harmonic 1-forms in dimension four from complex geometry}
\label{Sect_Z2HarmFormsFromCmplxGeom}

We begin by recalling a familiar picture of  a Riemann surface $S$ equipped with a branched double-cover $\pi:S \to \C$,
with branching set $\Delta:=\{ z_1,\dots, z_k,\dots \}$. The group $\Z_2$ acts freely on $S\setminus\pi^{-1}(\Delta)\to 
\C\setminus\Delta$, and $S\setminus\pi^{-1}(\Delta)$ is thereby a principal $\Z_2$ bundle. The associated real line bundle 
\[
\cI_\Delta :=\big (S\setminus\pi^{-1}(\Delta)\big )\times_{\Z_2}\R
\]
over $\C\setminus \Delta$ has monodromy $-1$ around each loop encircling $z_j$ and no other $z_i$. Denote by $\iota: S \to S$
the natural involution. A holomorphic function on $S$ which is odd with respect to $\iota$ is a holomorphic section $f$ of $\cI\otimes\C$. 
Notice the following:
\begin{itemize}[itemsep=-2pt,topsep=2pt]
	\item If $p\in\pi^{-1}(\Delta)$ is a fixed point of the $\Z_2$ action, then $f(\pi(p)) = 0$. 
	\item The real and imaginary  parts of $f$ are well-defined as sections of $\cI$.
\end{itemize}
The prototypical example is
\begin{equation*}
S =\{ (z,w)\in \C^2\mid w^2=z \} \qquad\text{with}\qquad \pi(z,w) = z,\ \ \iota(z,w) = (z,-w).
\end{equation*}
The holomorphic function $F(z,w) = w$ is identified with the holomorphic section $f(z)=\sqrt z$ of $\cI \otimes \C$;
$\Re \sqrt{z}$ and $\Im \sqrt{z}$ are real-valued harmonic sections of $\cI$.  

We generalize this example by first noting that if $X$ is a complex surface equipped with a holomorphic double branched cover 
$\pi: X \to \C^2$, with branching set an analytic subvariety $\Sigma\subset \C^2$, then there is a real Euclidean line bundle 
$\cI_\Sigma\to\C^2\setminus\Sigma$ with nontrivial monodromy around $\Sigma$ and a correspondence between odd 
holomorphic functions on $X$ and holomorphic sections of $\cI_\Sigma\otimes\C$. 

Our examples in dimension $4$ are all based on the
\begin{proposition} 
Suppose that $h\colon \C^2\to\C$ is any holomorphic function, which vanishes to odd order along a subvariety $\Sigma$. 
Then 
\begin{equation*}
X:=\big \{ (y, z,w)\in \C^3\mid y^2=h(z,w)^3 \big \}
\end{equation*}
is a K\"ahler surface which double covers $\C^2$ via the projection $(y,z,w) \mapsto (z,w)$. The ramification set is $\Sigma= h^{-1}(0)$, and the $\Z_2$
action on $X$ and real Euclidean line bundle $\cI$ over $\C^2 \setminus \Sigma$ are defined as above. The surface $X$
is the graph of the $2$-valued holomorphic map $h(z,w)^{3/2}$ branched over $\Sigma$.  Then $h(z,w)^{3/2}$ is a holomorphic 
section of $\cI\otimes\C$, and in particular 
\[
f:=\Re h^{3/2}
\] 
is a harmonic (and in fact, pluriharmonic) section of $\cI$ with vanishing order at least $3/2$ along $\Sigma$.

Finally, 
\[
\om = df
\]
is a $\Z_2$ harmonic $1$-form on $\C^2$ with branching set $\Sigma$, and with vanishing order at least $1/2$ there.
\end{proposition}

\begin{remark}
As noted, the power $3/2$ can be replaced by any number of the form $(2k+1)/2$, $k \in \Z$. However, the corresponding
function $|\om|$ is continuous on $\C^2$ only if $k \geq 1$. 
\end{remark}

In the same manner we obtain families of $\Z_2$ harmonic functions and $1$-forms.
\begin{corollary} 
Let $h=h(a, z, w)$ be a smooth function on $\mathcal U \times \C^2$, where $\mathcal U\subset \R^N$ is a parameter set, such that each $h_a:=h(a,\cdot,\cdot)$ is a holomorphic function on $\C^2$.  
The zero locus 
\begin{equation}
\label{Eq_SigmaA}
\Sigma_a=\{ (z,w) \mid h_a(z,w)= 0 \} \subset \C^2
\end{equation}
is smooth for generic values of $a$, e.g., when $\nabla h_a \neq 0$. Suppose, as above, that $h_a$ vanishes
to odd order on $\Sigma_a$. Denote by $\cI_a$ the associated real Euclidean line bundle $\cI_a$ with 
nontrivial monodromy around $\Sigma_a$. Then
\begin{equation}
\label{Eq_Basic2valHolFn_SmoothDeform}
F_a(z,w):= h(a,z,w)^{3/2},\qquad \ f_a = \Re h(a,z,w)^{3/2},\qquad \mbox{and}\quad  \om_a = d f_a
\end{equation}
are holomorphic sections of $\cI_a\otimes \C$ and harmonic sections and $1$-forms valued in $\cI_a$, respectively. 
\end{corollary}

\begin{remark}
We consider in particular two types of families: first, the family $h(a,z,w) = h(z,w) - a$, $a \in \C$,
and second, $h(a,z,w) = h \circ \varphi_a (z,w)$ where $\varphi_a$ is any family
of unitary transformations of $\R^4 \cong \C^2$.  
Notice that we can allow $\varphi_a$ to be a family of isometries only: while in this case $F_a$ does not need to be holomorphic, $f_a$ and $\om_a$ remain harmonic. 
\end{remark}

Based on these general claims, we now list several examples and describe the new phenomena exhibited in each.

\begin{example} (Desingularizing an immersion)
\label{Ex_UnionOfCoordAxes}
The function $h_0(z,w) = zw$ has zero locus $\Sigma_0$ the union of the two coordinate axes. Thus (omitting the holomorphic intermediary),
\begin{equation}
\label{Eq_Basic2valHolFn}
\Re h_0(z,w)^{3/2} =\Re \big ( z^{3/2} w^{3/2} \big )\ \mbox{and}\ \ \om_0 = d \Re h_0(z,w)^{3/2} = \frac32 \Re \Big(
h^{\frac12} \Big (\frac{\del h}{\del z} \, dz + \frac{\del h}{\del w}\, dw\Big) \Big)
\end{equation}
are $\Z_2$ harmonic functions and $1$-forms branching around $\Sigma_0$. 

Now consider the family of functions $h_{a,b,c}(z,w) = (z-b)(z-c) - a$. The zero locus $\Sigma_{a,b,c}$ is the smooth
curve $(z-b)(w-c) = a$ when $a \neq 0$ and the union of complex lines $\{z=b\} \cup \{w=c\}$ when $a=0$. We obtain
the family of $\Z_2$ harmonic functions and $1$-forms
\[
f_{a,b,c}(z,w) = \Re h_{a,b,c}(z,w)^{3/2}, \ \ \mbox{and}\ \ \ \om_{a,b,c} = d \Re h_{a,b,c}(z,w)^{3/2}.
\]

This family has branching locus a smooth curve when the parameters lie in the Zariski open set $\{a \neq 0\}$,
and are nodal curves when $a = 0$. 
\end{example}



Returning to the more general formulation of the Proposition, we obtain $\Z_2$ harmonic functions and $1$-forms
which branch along any complex analytic subvariety of (complex) codimension $1$.  These sets can have many different
types of singular behavior; thus (at least in dimension four), singularities of branching sets of $\Z_2$ harmonic
$1$-forms are at least as complicated as those of one dimensional complex analytic sets.   Note that for all of these
examples in four dimensions, the singularities of the branching sets $\Sigma_a$ are isolated points.



\begin{example}[$\Sigma$ is any collection of complex lines through the origin]
\label{Ex_TangCones}
Let 
\begin{equation*}
\Sigma=\bigcup_{j=1}^J\Sigma_j,\quad\text{where}\ \  \Sigma_j=\{ a_j z + b_j w=0 \},\ (a_j, b_j) \neq (0,0)\ \forall\, j
\end{equation*}
be any finite collection of complex lines through the origin in $\C^2$. This is the zero locus of $h(z,w)= \prod (a_j z + b_j w)$. 
Hence by the Proposition above, there is a $\Z_2$ harmonic $1$-form with branching locus $\Sigma$.  The unique singular
point of $\Sigma$ is the origin, and the tangent cone of $\Sigma$ at that point is $\Sigma$ itself.  This shows that
in real dimension four, the set of all possible tangent cones of branching sets of  $\Z_2$ harmonic 1-form has positive dimension 
(even modulo isometries).  We remark that the corresponding statement is unknown in dimension $3$, but may not be true.
The examples constructed in \cite{TaubesWu20_Examples}, where the branching locus is a union of rays emanating from
the origin in $\R^3$, appears to be rigid (though this is not proved). 
\end{example}
\begin{example}(Ramified branched covers) Set $h_a(z,w):=w^2-a(z^3 + 1)$, $a\in \C$. When $a\neq 0$, $\Sigma_a:=h_a^{-1}(0)$ is 
an (affine) elliptic curve which is a ramified double covering of the plane $w=0$. When $a=0$, $\Sigma_0$ is the plane $w=0$ itself. 

Families of branching sets can thus converge to limits which have multiplicity greater than $1$.
\end{example}

\section{Examples in dimension $3$}
We now discuss a number of different examples and constructions of $\Z_2$ harmonic $1$-forms in $3$-manifolds; some of these
were known previously, but are included here to round out this compendium.  The first two subsections describe constructions which
still take some advantage of complex geometry. The other examples are of a somewhat different nature.

\subsection{A $\Z_2$ harmonic form with a non-constant vanishing order along the branching set}

Let $h$ be a $\Z_2$ harmonic function on $\R^2\setminus\{ 0\}$ with values in the unique non-trivial real line bundle. If $(x,y,z)$ are standard coordinates on $\R^3$, $H(x,y,z):= z\,h(x,y)$ is a $\Z_2$ harmonic function on $\R^3\setminus \{z\text{-axis}\}$. For example, denoting $w:=x+iy$ and choosing $h(w)=\Re (w^{3/2})$, we obtain $H=z\Re (w^{3/2})$. In particular, 
\begin{equation}
	\label{Eq_Z2HR3}
	\om:= 2 \Re\big (  w^{\frac 32}  \big )\, dz + 3z\Re \big ( w^{\frac 12}\, dw\big )
\end{equation}
is a $\Z_2$ harmonic 1-form on $\R^3\setminus \{z\text{-axis}\}$. 
This has the interesting feature that $\om$ vanishes to order $\frac 32$ at the origin, but only to order $\frac 12$ at every 
other point of the $z$-axis.  To our knowledge, this the only known example with this property.   This example shows that
the nondegeneracy hypothesis in the main theorems of ~\cite{Takahashi15_Z2HarmSpinors_Arx}, \cite{Parker23_DeformZ2HarmSpinors_Arx}
and \cite{DonaldsonQJM:21} is a genuine restriction.

\subsection{Harmonic morphisms}
We have already noted the (trivial) fact that if $\om$ is a $\Z_2$ harmonic $1$-form, then its pullback with respect to any isometry
of the ambient space is also $\Z_2$ harmonic.  There is an interesting sort of generalization of this to any dimension. 

Let $M$ and $N$ be two Riemannian manifolds. Recall from \cite{BairdWood03_HarmMorphisms} that a map $\varphi\colon M\to N$ 
is called a \emph{harmonic morphism} if, for any locally defined harmonic function $h$ on $N$, the pullback $\varphi^*h=h\comp \varphi$ 
is harmonic on $M$.  It follows directly from this that if $\zeta$ is a locally defined harmonic $1$-form on $N$, then its pullback
$\varphi^*\zeta$ is harmonic on $M$. Still slightly more generally, if $(\om_N, \Sigma_N, \cI_N)$ is a $\Z_2$ harmonic $1$-form, 
then $\big (\varphi^*\om_N,\; \varphi^{-1}(\Sigma_N), \; \varphi^*\cI_N\big )$ is a $\Z_2$ harmonic 1-form on $M$. 

Unfortunately the condition of being a harmonic morphism is quite rigid and there are relatively few examples, see
\cite{BairdWood03_HarmMorphisms} for further details. However, there are some nontrivial ones which can be 
used to construct new $\Z_2$ harmonic $1$-forms in $3$ dimensions.

\begin{example}
If both $M$ and $N$ are K\"ahler, then any holomorphic map $\varphi$ is a harmonic morphism, see \cite{BairdWood03_HarmMorphisms}*{Cor.\,8.1.6}. 
In particular, given a holomorphic function $h\colon \C^2\to\C$, the pull-back of the elementary $\Z_2$ harmonic $1$-form $\Re (z^{1/2}dz)$ 
on $\C$ is simply $\frac23 d (\Re h(w_1, w_2)^{3/2})$, the main example of Section 2. 
\end{example}

\begin{example}[Linked $\Sigma$]
\label{Ex_HopfMap}
The Hopf map
\begin{equation*}
\varphi\colon S^3\to S^2,\qquad \varphi(z,w) = \big (  |z|^2-|w|^2,\; 2z\bar w \big )
\end{equation*} 
is a harmonic morphism, cf.\ \cite{BairdWood03_HarmMorphisms}*{Ex.\,5.6.4}.  Here $S^3 \subset \C^2$ and $S^2 \subset \R \oplus \C$
have their standard round metrics. 
The preimage by $\varphi$ of any point in $S^2$ is a great circle in $S^3$, and the preimage of any pair of points in $S^2$ is a Hopf link. 
	
Identify $S^2\setminus\{ \infty \}$ with $\C$; for any holomorphic function $p(z)$ on $\C$, set
\begin{equation}
\label{Eq_LocZ2HarmFormC}
\om:= \Re\big (  p(z)^{\frac 12}\, dz  \big ).
\end{equation} 
Thus $\om$ is a $\Z_2$ harmonic $1$-form on $\C \setminus p^{-1}(0)$, where the twisting bundle $\cI_p$ has nontrivial monodromy 
around each zero of $p$ of odd multiplicity. To simplify the notations, assume every zero of $p$ has odd multiplicity. Then
\begin{equation*}
\Big (  \varphi^*\om,\; \varphi^{-1}\big (p^{-1}(0)\big ),\; \varphi^*\cI_p) \Big )	
\end{equation*}
is a $\Z_2$ harmonic $1$-form on $S^3\setminus \varphi^{-1}(\infty)=S^3\setminus S^1$ with its round metric. The behavior of $\varphi^* \om$
at this omitted circle depends on the behavior of $p(z)$ at infinity, so we regard this only as a `local' example, on $S^3$ away from this circle.
The branching set of this harmonic $1$-form is $\varphi^{-1}\big (p^{-1}(0)\big )$, which is a finite collection of great circles, any pair of which
is a Hopf link.   These may be the first examples of (local) $\Z_2$ harmonic 1-forms where $\Sigma$ is a nontrivial link.
\end{example}

\begin{example}[Knotted branching sets and multiple covers]
\label{Ex_TorusKnots}
Our next examples relies on a family of harmonic morphisms $S^3\to S^2$ generalizing the Hopf map.  Regard $S^3 \subset \C^2$ as usual.
Fix a pair $p$ and $q$ of coprime integers, and define
\begin{equation*}
\pi_{p,q}\colon S^3\to \CP^1\cong S^2; \quad \pi_{p,q}\big (  z_1, z_2\big ) = \big [ \, z_1^p: z_2^q\, \big].
\end{equation*} 
This is the standard Hopf map when $p=q=1$. 

The map $\pi_{p,q}: S^3 \to S^2$ is a Seifert fibration with two singular fibers, the preimages of the north and south poles. The singular fibers 
are unknots, whereas any other fiber is a torus knot, of type $(p,q)$ around one singular fiber and of type $(q,p)$ around the other.

Let $g_{p,q}$ be the restriction of the metric $p^2|dz_1|^2 + q^2|dz_2|^2$ to $S^3$. The resulting Riemannian manifold $Q_{p,q}$ is an ellipsoid.
By \cite[Ex.\,10.4.2]{BairdWood03_HarmMorphisms}, there is a diffeomorphism $f\colon S^2\to S^2$ such that $\varphi_{p,q}:=
f\comp \pi_{p,q}\colon Q_{p,q}\to S^2$ is a  harmonic morphism.

For any $a\in \C$, set  $\om_a:=\Re\big ( (z-a)^{1/2}dz \big )$. Choose the embedding $\C\subset S^2$ so that the origin in $\C$ 
corresponds to $f(-1,0)$, where $(-1,0)\in S^2\subset \R\oplus\C$. Then if $a\in  
\C\setminus\{ 0 \}$, the triple 
$\big (\varphi_{p,q}^*\,\om_a , \varphi_{p,q}^{-1}(a), \varphi_{p,q}^*\,\cI_a \big )$ is a (local) $\Z_2$ harmonic 1-form on $Q_{p,q}$ with
$\Sigma_a=\varphi^{-1}(a)$ a $(p,q)$ torus knot (here, $\Sigma_a$ passes $q$ times along the core $\Sigma_0$).
	
This example has the following interesting property: if $a \approx 0$, then  $\Sigma_a$ is close to the unknot $\Sigma_0$. However,
$\Sigma_a$ is topologically a $q$-fold covering of $\Sigma_0$. In other words, $\big (\varphi_{p,q}^*\,\om_a , \varphi_{p,q}^{-1}(a), 
\varphi_{p,q}^*\,\cI_a \big )$ is a family of $\Z_2$ harmonic $1$-forms such that for each $a\neq 0$ the corresponding branching set 
$\Sigma_a$ is a $q$-fold covering of $\Sigma_0$.  
\end{example}

As in Example~\ref{Ex_HopfMap}, this example is local since the $\Z_2$ harmonic $1$-form constructed here typically blows up along 
one of the fibres.  Nevertheless, these two examples are explicit, and provide a vivid visualization of this phenomena.

The constructions of Examples~\ref{Ex_HopfMap} and~\ref{Ex_TorusKnots} can in fact be modified to obtain global examples. 
We first outline a general approach and then present some concrete examples below. 

Let $\pi\colon M^3\to F$ be a Seifert fibration, where $F$ is  a smooth oriented surface. 
Choose a Riemannian metric, and hence a complex structure, on $F$. Then $M$ admits a Riemannian metric $g$ such that 
$\pi$ is a harmonic morphism~\cite{Baird90_HarmMorphismsAndS1Actions}*{P.\,196}. Furthermore, given a holomorphic 
quadratic differential $q$ on $F$, we can construct a $\Z_2$ harmonic 1-form on $F$ by setting $\om_q:=\Re \sqrt q$, 
see~\cite{Taubes13_PSL2Ccmpt}. The branching set of $\om_q$ consists of the union of odd multiplicity zeros of $q$. 
Consequently, $\pi^*\om_q$ is a $\Z_2$ harmonic $1$-form on $M$ with branching set $\Sigma_q=\pi^{-1}\big (q_{odd}^{-1}(0)\big )$.
Any surface $F$ of genus at least $2$ admits holomorphic quadratic differentials with nonempty zero set.  We can therefore construct 
global $\Z_2$ harmonic $1$-forms on $3$-manifolds.  

Here is an illustration: 
\begin{example}
The Brieskorn manifold $M (a_1, a_2, a_3)$ is the link of the isolated singularity in the complex surface 
$V(a_1, a_2, a_3):=\{ z_1^{a_1} + z_2^{a_2} + z_3^{a_3} =0 \}\subset \C^3$, i.e., 
\begin{equation}
M (a_1, a_2, a_3)=V(a_1, a_2, a_3)\cap S^5. 
\end{equation}
For simplicity we restrict to the case $a_1=a_2=a_3=:a\ge 4$. Then the standard $\U(1)$-action on $\C^3$ preserves $M (a, a, a)$ 
and the quotient is the surface 
\begin{equation*}
F:=\big \{ [z_1:z_2:z_3]\in \CP^2\mid z_1^a + z_2^a + z_3^a =0  \big \};
\end{equation*} 
this has genus $(a-1)(a-2)/2\ge 3$.

It is known that the fibers of the natural projection $\pi\colon M(a,a,a)\to F$ are geodesics with respect to the inherited metric 
on $M(a,a,a)$, and furthermore, $\pi$ is horizontally weakly conformal~\cite{Baird90_HarmMorphismsAndS1Actions}*{Def.\,2.4.2}.
Hence by~\cite{Baird90_HarmMorphismsAndS1Actions}*{Thm.\,4.5.4}, $\pi$ is a harmonic morphism. Topologically, $M(a,a,a)$ is 
a circle bundle over $F$ with Euler class $-a$.

If $q$ is a nontrivial holomorphic quadratic differential $q$ on $F$, then $\pi^*\om_q$ is a $\Z_2$ harmonic $1$-form on $M(a,a,a)$ with 
branching set $\Sigma = \pi^{-1}\big (q_{odd}^{-1}(0)\big )$, which is a link in $M(a,a,a)$.  
\end{example}

\subsection{Smoothing immersions}
\begin{figure}[ht]
\centering
\includegraphics[width=0.7\linewidth]{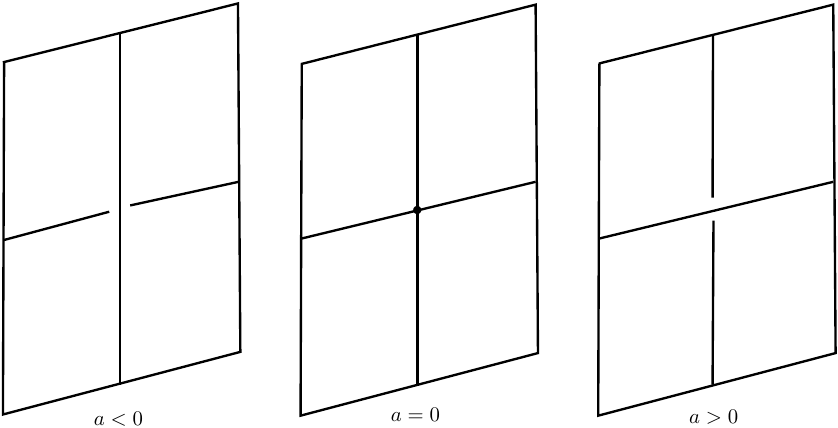}
\caption{Branching locus changes from undercrossing to overcrossing}
\label{FigUnderOverCross}
\end{figure}

\begin{figure}[ht]
	\centering
	\includegraphics[width=0.7\linewidth]{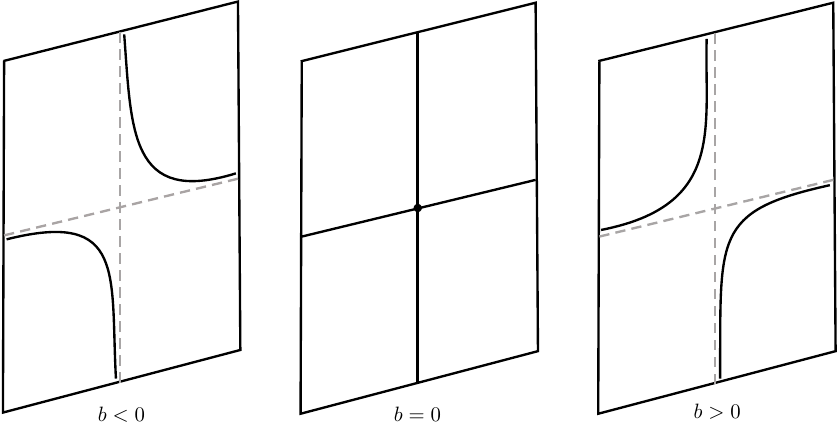}
	\caption{A smoothing of the branching set}
	\label{FigSmoothing}
\end{figure}

We have already exhibited families of $\Z_2$ harmonic $1$-forms $\om_t$ on $\C^2$ where the branching set $\Sigma_0$
is immersed and the nearby $\Sigma_t$ are embedded.  It is reasonable to expect similar phenomena in three dimensions, but
there are no known examples yet. We conjecture, however, that there are $\Z_2$ harmonic $1$-forms $(\om, \Sigma, \cI)$   on $\R^3$ such that:
\begin{itemize}
\item $\Sigma$ is the union of two intersecting lines;
\item There exist two smooth families of $\Z_2$ harmonic $1$-forms $(\om_t, \Sigma_t, \cI_t)$ and $(\hat \om_t, \hat\Sigma_t, \hat \cI_t)$, 
$|t| < \epsilon$, such that: 
\begin{itemize}
\item $(\om_t, \Sigma_t, \cI_t)|_{t=0} =(\hat \om_t, \hat \Sigma_t, \hat \cI_t)|_{t=0} = (\om, \Sigma, \cI)$;
\item If $t\neq 0$, $\Sigma_t$ is either an under- or over-crossing as in Figure\, \ref{FigUnderOverCross};
\item If $t\neq 0$, then $\hat \Sigma_t$ is a smoothing as in Figure\,\ref{FigSmoothing}.
\end{itemize}
\end{itemize}
Notice that over- and under-crossings are well-defined only when  $\Sigma_t$ is endowed with an orientation.

If such a family of $\Z_2$ harmonic $1$-forms $\om_t$ were to exist, where $\Sigma_t$ is a smooth desingularization of $\Sigma_0$,
then it can be used as one of the building blocks in a gluing construction to obtain families of $\Z_2$ harmonic $1$-forms with
branching sets $\Sigma_t$ smoothings of more general immersed curves $\Sigma_0$.  This is being carried out in current ongoing
work by the second author and Siqi He. As in other problems of this type, there is a loss of regularity which necessitates
the use of the Nash implicit function theorem.

\subsection{$\Z_2$-harmonic $1$-forms with symmetry}
We conclude this section by recalling a clever construction by Sun \cite[Appendix B]{Sun22} which produces $\Z_2$ harmonic functions
and $1$-forms on $\R^n$ with (possibly large) polynomial growth at infinity, and which have a certain rotational symmetry.   

The first step is a more general one which starts with any compact smooth submanifold $\Sigma$ of dimension $n-2$ in $\R^n$ and
a harmonic polynomial $p(x)$ and produces a harmonic section of the twisted real Euclidean line bundle $\cI$ with branching set $\Sigma$.  
Sections constructed this way typically only decay like the square root of distance to $\Sigma$, hence their differentials blow up to 
order $-1/2$ at $\Sigma$.  These are unsuitable for our goals. However, in the second step we show how to use additional 
symmetry hypotheses to find examples of such harmonic sections which decay to the desired order $3/2$ at $\Sigma$. 

Sun's paper \cite{Sun22} contains a somewhat lengthy proof of the first step, but notes that the second author here pointed
out to him a more standard general procedure. We record this alternate proof here, but emphasize that the idea for the 
existence of such harmonic sections is entirely due to Sun.

Start then with the submanifold $\Sigma \subset \R^n$, $n\ge 3$, and the associated real line bundle $\cI$. Consider the singular
manifold $X$ which is obtained by gluing two copies of $\R^n \setminus \Sigma$ along $\Sigma$ so that the familiar
branching occurs in directions normal to $\Sigma$ and $X$ is connected. Thus $X \setminus \Sigma$ is identified with
the set of all unit norm sections of $\cI$, and $X$ itself as the multi-valued graph of a section $\sigma$ of $\cI$ which 
equals $\pm r^{1/2}$ near $\Sigma$, with $\sigma = \pm 1$ outside a somewhat larger neighborhood, and vanishes nowhere 
away from $\Sigma$. Denote by $\pi\colon X \to \R^n$ the natural projection. If the constant $R_0$ is sufficiently large,
 then $\pi^{-1}( \R^n \setminus B_{R_0}(0))$ consists of two sheets $X^+_{R_0} \sqcup X^-_{R_0}$.

Next, let $p(x)$ be any harmonic polynomial on $\R^n$. Choose a 
smooth cutoff function $\chi(x)$ on $\R^n$ which equals $1$ outside some large ball $B_{R_2}(0)$, and which vanishes 
inside a somewhat smaller ball $B_{R_1}(0)$ which contains $\Sigma$. Consider the function $U$ on $X$ which equals 
$\pm \chi(x) p(x)$ on the two sheets $X^\pm_{R_1}$ and which vanishes on the remaining part of $X$. Define $H = \Delta U$; this
is smooth and compactly supported.

We now find a solution to the equation $\Delta V = H$ globally on $X$ where $V$ grows at infinity at a slower rate than $U$,
and decays at least like $\mathrm{dist}_\Sigma^{1/2}$ near $\Sigma$.  If we can do so, then $u = U - V$ is the harmonic 
function on $X$ we are seeking in this first step. 

To solve this problem, we define suitable function spaces. All of this is standard material but included for the reader's 
convenience.  First define the weighted $L^2$ spaces $M^\delta L^2$, where $M$ is a smooth strictly positive 
function which equals $1$ on $\pi^*(B_{R_0}(0))$ and $\pi^* (1+|z|^2)^{1/2}$ outside a slightly larger set $\pi^*(B_{R_1}(0))$.
Denoting $\hat \chi:= \pi^*\chi$, define the space
\[
H^2_{e,\delta}(X) = \big \{u \mid \nabla^j (\hat\chi\, u )\in M^{\delta-j} L^2(X),\ 0 \leq j \leq 2, \text{ and } (1-\hat\chi)u\in H^2_e(X) \big\},
\]
where $H^2_e(X)$ denotes the edge Sobolev space. We refer the reader to~\cite{Mazzeo91_EllTheoryOfDiffEdgeOp} for the definition of the edge Sobolev spaces and for further details of edge calculus. We just note that near $\Sigma$ one can think of $X$ along each normal fiber as a cone of angle $4\pi$ and the space $H^2_{e}$ is adapted to this situation. 

Consider the space
\[
H^2_{\delta, D} (X) = \big \{u \in H^2_{e,\delta}\mid\  \nabla \big ((1-\hat\chi ) u \big) \in L^2(X) \big \}
\]
This imposes the Fridrichs extension condition near $\Sigma$ and corresponds to standard weighted Sobolev spaces near $\infty$. 
This domain admits functions decaying like $r^{1/2}$ but exclude 
those which blow up like $r^{-1/2}$ near $\Sigma$. 


%

Now consider the map
\[
\Delta: H^2_{\delta, D}(X) \longrightarrow   M^{\delta-2} L^2(X).
\]
It is well-known that this map is Fredholm for every $\delta \not\in \Lambda := \{ 2-\frac{n}{2} - k, \frac{n}{2} + k \mid k = 0, 1, 2, \ldots\}$; 
furthermore, this map is injective when $\delta < n/2$ and surjective when $\delta > 2 - n/2$.   To check this, we note that if $\Delta u = 0$
and $u \in M^\delta L^2$ with $\delta < n/2$, then by a standard spherical harmonic expansion, $|u| \leq C M^{2-n}$ near infinity,
and thus $0 = \int u \Delta u = -\int |\nabla u|^2$, so $u \equiv 0$. The integration by parts is justified by the decay near infinity and
the Friedrichs condition near $\Sigma$.  The surjectivity statement for $\delta > 2-n/2$ then follows by a duality argument. 

Using this result, we then fix any $\delta > 2-n/2$ and then find a solution $V \in H^2_{\delta,D}$ to $\Delta V = H$.  So long 
as $\delta < n/2 + \mathrm{deg}\, p$, the harmonic section $u = U - V$ is asymptotic to $U$ as $|x| \to \infty$. It also decays
only at rate $r^{1/2}$ at $\Sigma$. 

Now let us turn to the second step, specializing to $n=3$ since this case is particularly transparent. Using coordinates
$(x_1,x_2,x_3)$ on $\R^3$, define $\Sigma = \{x_1^2 + x_2^2 = 1,\ x_3 = 0\}$. We now work entirely within the space of 
functions and sections which are invariant under rotations around the $x_3$ axis; this is consistent with the rotation symmetry
of $\Sigma$.  Performing the construction from step 1 using only harmonic polynomials $p(x)$ which enjoy a similar
rotation invariance, then the resulting harmonic functions $u(x) = u_p(x)$ are also invariant with respect to these same rotations.
There is a unique homogeneous harmonic polynomial of degree $k$ with this symmetry. It is called the $k^{\mathrm{th}}$
zonal harmonic, and is explicitly given as $\rho^k P_k(\cos \phi)$, where $\rho$ is the polar radial variable in $\R^3$,
$\phi$ is the azimuthal angle in spherical coordinates, and $P_k$ is the Legendre polynomial of degree $k$. 
Denote its ($1$-dimensional) span by $\mathcal P_k$.

Finally, as shown in essentially any of the papers cited earlier about $\Z_2$ harmonic functions, the harmonic section $u$ has
an expansion near $\Sigma$ in powers of $r = \dist(\Sigma, \cdot)$ and the angular coordinate $\theta$ in the fibers of the 
normal bundle around $\Sigma$ of the form
\[
u_p(x) = A_1^+(p) \cos \theta/2 \, r^{1/2} + A_1^-(p) \sin \theta/2\, r^{1/2} + \mathcal O(r^{3/2}).
\]
The leading coefficients are written as $A_1^\pm(p)$ to emphasize their dependence on $p$, and it is clear that this dependence is linear.
If we had not imposed the rotation invariance, these coefficients would also depend on the variable $t$ along $\Sigma$, but with this
rotation invariance these are simply constant in $t$.   We have thus constructed two linear functionals 
\[
\mathcal P  = \oplus_{k \geq 0} \mathcal P^k \longrightarrow \R^2, 
\quad p(x) \longmapsto (A_1^+(p), A_1^-(p)).
\]
Since $\dim \mathcal P = \infty$, the nullspace of this map is also infinite dimensional.  Any harmonic polynomial in this nullspace
yields, by the construction above, a $\Z_2$ harmonic function on $X$ which decays at least like $r^{3/2}$ at $\Sigma$;
the corresponding differentials are $\Z_2$ harmonic $1$-forms which branch along this $\Sigma$. 

As before, we attribute this observation about how to use rotational symmetry to Sun~\cite{Sun22}.

\section{Conclusion}  We hope that this collection of constructions and
examples might stimulate the reader toward this topic.  Intriguingly, it
appears to be much easier to construct this variety of examples in the $\Z_2$ 
harmonic $1$-form case, but not for $\Z_2$ harmonic spinors.  This is 
a good challenge. We conclude by mentioning a forthcoming paper by
He and Parker \cite{He-Parker}, in which a gluing theorem is proved for
$\Z_2$ harmonic spinors and $1$-forms. We expect this to be a useful
tool in understanding the landscape of these objects better. We also mention
the analytic deformation theory  for $\Z_2$ harmonic spinors developed by
the third author, and more recently by Parker, \cite{Takahashi15_Z2HarmSpinors_Arx, Parker23_DeformZ2HarmSpinors_Arx}. These theorems only cover the case where 
the branching set is a smooth curve in a $3$-manifold.  There is interest and considerable hope for being able to obtain a similar deformation theory when
the branching sets are immersed curves, or more generally, smoothly embedded 
graphs.  At the heart of these nonlinear deformation theorems are index theorems
for the linearized deformation operators. The index theorem in the smooth
branching case was obtained by Takahashi, and in the case where the branching
set is a smoothly embedded graph, the index theorem has been obtained recently by the three authors of this note~\cite{HMT23-IndexGraphs-Arx}.

\bibliography{references}

\def\cprime{$'$}
\begin{bibdiv}
\begin{biblist}

\bib{Baird90_HarmMorphismsAndS1Actions}{article}{
      author={Baird, P.},
       title={Harmonic morphisms and circle actions on {$3$}- and
  {$4$}-manifolds},
        date={1990},
        ISSN={0373-0956},
     journal={Ann. Inst. Fourier (Grenoble)},
      volume={40},
      number={1},
       pages={177\ndash 212},
         url={http://www.numdam.org/item?id=AIF_1990__40_1_177_0},
}

\bib{BairdWood03_HarmMorphisms}{book}{
      author={Baird, P.},
      author={Wood, J.},
       title={Harmonic morphisms between {R}iemannian manifolds},
      series={London Mathematical Society Monographs. New Series},
   publisher={The Clarendon Press, Oxford University Press, Oxford},
        date={2003},
      volume={29},
        ISBN={0-19-850362-8},
         url={https://doi.org/10.1093/acprof:oso/9780198503620.001.0001},
      review={\MR{2044031}},
}

\bib{Donaldson17_AdiabLimitsG2}{incollection}{
      author={Donaldson, S.},
       title={Adiabatic limits of co-associative {K}ovalev--{L}efschetz
  fibrations},
        date={2017},
   booktitle={Algebra, geometry, and physics in the 21st century},
      series={Progr. Math.},
      volume={324},
   publisher={Birkh\"{a}user/Springer, Cham},
       pages={1\ndash 29},
         url={https://doi.org/10.1007/978-3-319-59939-7_1},
      review={\MR{3702382}},
}

\bib{DonaldsonQJM:21}{article}{
      author={Donaldson, S.K.},
       title={{Deformations of multivalued harmonic functions}},
        date={2021},
     journal={Quarterly J. Math},
      volume={72},
      number={1-2},
       pages={199\ndash 235},
}

\bib{Haydys22_SWPSL2R}{article}{
      author={Haydys, A.},
       title={Seiberg--{W}itten monopoles and flat {${\rm PSL}(2,\mathbb
  R)$}-connections},
        date={2022},
        ISSN={0001-8708,1090-2082},
     journal={Adv. Math.},
      volume={409},
       pages={Paper No. 108686, 18},
         url={https://doi.org/10.1016/j.aim.2022.108686},
      review={\MR{4483244}},
}

\bib{HMT23-IndexGraphs-Arx}{article}{
      author={Haydys, A.},
      author={Mazzeo, R.},
      author={Takahashi, R.},
       title={An index theorem for {Z}/2-harmonic spinors branching along a
  graph},
        date={2023},
      eprint={\arxiv{2310.15295}},
}

\bib{HaydysWalpuski15_CompThm_GAFA}{article}{
      author={Haydys, A.},
      author={Walpuski, T.},
       title={A compactness theorem for the {S}eiberg--{W}itten equation with
  multiple spinors in dimension three},
        date={2015},
        ISSN={1016-443X},
     journal={Geom. Funct. Anal.},
      volume={25},
      number={6},
       pages={1799\ndash 1821},
         url={http://dx.doi.org/10.1007/s00039-015-0346-3},
        note={Erratum available via
  \href{http://dx.doi.org/10.13140/RG.2.1.2559.1285 }{doi:
  10.13140/RG.2.1.2559.1285}},
      review={\MR{3432158}},
}

\bib{He23_BranchedDeform}{article}{
      author={He, S.},
       title={The branched deformations of the special lagrangian
  submanifolds},
        date={2023},
     journal={Geom. Funct. Anal.},
      volume={25},
         url={https://doi.org/10.1007/s00039-023-00645-8},
}

\bib{He-Parker}{article}{
      author={He, S.},
      author={Parker, G.},
       title={Connect sum and torus sum formulas for $\mathbb{Z}_2$ harmonic
  spinors and $1$-forms on $3$-manifolds},
        date={2024},
        note={To Appear},
}

\bib{Mazzeo91_EllTheoryOfDiffEdgeOp}{article}{
      author={Mazzeo, R.},
       title={Elliptic theory of differential edge operators. {I}},
        date={1991},
        ISSN={0360-5302},
     journal={Comm. Partial Differential Equations},
      volume={16},
      number={10},
       pages={1615\ndash 1664},
         url={https://doi.org/10.1080/03605309108820815},
      review={\MR{1133743}},
}

\bib{Parker23_DeformZ2HarmSpinors_Arx}{article}{
      author={Parker, G.},
       title={Deformations of $\mathbb {Z}_2$-harmonic spinors on 3-manifolds},
        date={2023},
      eprint={\arxiv{2301.06245}},
}

\bib{Sun22}{article}{
      author={Sun, W.},
       title={On $\mathbb {Z}_2$ harmonic functions on $\mathbb {R}^2$ and the
  polynomial {P}ell's equation},
        date={2022},
      eprint={\arxiv{2209.11893}},
}

\bib{Takahashi15_Z2HarmSpinors_Arx}{article}{
      author={Takahashi, R.},
       title={The moduli space of {$S^1$}-type zero loci for {Z}/2-harmonic
  spinors in dimension 3},
        date={2015},
      eprint={\arxiv{1503.00767}},
}

\bib{Taubes13_CxASD_Arx}{article}{
      author={Taubes, C.},
       title={Compactness theorems for {$\mathrm{SL}(2;{\Bbb C})$}
  generalizations of the 4-dimensional anti-self dual equations},
        date={2013},
      eprint={\arxiv{1307.6447}},
}

\bib{Taubes13_PSL2Ccmpt}{article}{
      author={Taubes, C.},
       title={{${\rm PSL}(2;\Bbb C)$} connections on 3-manifolds with {${\rm
  L}^2$} bounds on curvature},
        date={2013},
        ISSN={2168-0930},
     journal={Camb. J. Math.},
      volume={1},
      number={2},
       pages={239\ndash 397},
         url={http://dx.doi.org/10.4310/CJM.2013.v1.n2.a2},
        note={Corrigendum: Camb. J. Math. \textbf{3} (2015), no. 4, 619--631.},
      review={\MR{3272050}},
}

\bib{Taubes14_ZeroLoci_Arx}{article}{
      author={Taubes, C.},
       title={The zero loci of {Z}/2 harmonic spinors in dimension 2, 3 and 4},
        date={2014},
      eprint={\arxiv{1407.6206}},
}

\bib{Taubes15_CorrigendumToPSL2C}{article}{
      author={Taubes, C.},
       title={Corrigendum to ``{$PSL(2;\Bbb C)$} connections on 3-manifolds
  with {$L^2$} bounds on curvature'' [{C}ambridge {J}ournal of {M}athematics
  1(2013) 239--397] [ {MR}3272050]},
        date={2015},
        ISSN={2168-0930},
     journal={Camb. J. Math.},
      volume={3},
      number={4},
       pages={619\ndash 631},
      review={\MR{3435274}},
}

\bib{Taubes16_SWDim4_Arx}{article}{
      author={Taubes, C.},
       title={On the behavior of sequences of solutions to {$\rm U(1)$}
  {S}eiberg--{W}itten systems in dimension 4},
        date={2016},
      eprint={\arxiv{1610.07163}},
}

\bib{Taubes17_VafaWitten_Arx}{article}{
      author={Taubes, C.},
       title={The behavior of sequences of solutions to the {V}afa--{W}itten
  equations},
        date={2017},
      eprint={\arxiv{1702.04610}},
}

\bib{Taubes18_SequencesKapustWitten_Arx}{article}{
      author={Taubes, C.},
       title={Sequences of nahm pole solutions to the $\mathrm{SU}(2)$
  {K}apustin--{W}itten equations},
        date={2018},
      eprint={\arxiv{1805.02773}},
}

\bib{TaubesWu20_Examples}{inproceedings}{
      author={Taubes, C.~H.},
      author={Wu, Y.},
       title={Examples of singularity models for {$\mathbb Z/2$} harmonic
  1-forms and spinors in dimensional three},
        date={2020},
   booktitle={Proceedings of the {G}\"{o}kova {G}eometry-{T}opology
  {C}onferences 2018/2019},
   publisher={Int. Press, Somerville, MA},
       pages={37\ndash 66},
      review={\MR{4251085}},
}

\bib{Zhang22_Rectifiability}{article}{
      author={Zhang, B.},
       title={Rectifiability and {M}inkowski bounds for the zero loci of
  {$\mathbb Z/2$} harmonic spinors in dimension 4},
        date={2022},
        ISSN={1019-8385,1944-9992},
     journal={Comm. Anal. Geom.},
      volume={30},
      number={7},
       pages={1633\ndash 1681},
      review={\MR{4596625}},
}

\end{biblist}
\end{bibdiv}

\end{document}